\documentclass[12pt,a4paper]{amsart}

\usepackage{amsmath}
\usepackage{amssymb}
\usepackage{latexsym}
\usepackage{xypic}
\input xy
\xyoption{all}

\newcommand{\dual}{\makebox[0mm]{}^{{\scriptstyle\vee}}}
\newtheorem{theorem}{Theorem}[section]
\newtheorem{lemma}[theorem]{Lemma}
\newtheorem{proposition}[theorem]{Proposition}
\newtheorem{definition}[theorem]{Definition}
\newtheorem{corollary}[theorem]{Corollary}

\newtheorem{exmp}[theorem]{Example}
\newtheorem{exmps}[theorem]{Examples}
\newtheorem{rema}[theorem]{Remark}
\newenvironment{rem}{\begin{rema}\rm}{\end{rema}}
\newenvironment{example}{\begin{exmp}\rm}{\end{exmp}}
\newenvironment{examples}{\begin{exmps}\rm}{\end{exmps}}
\newenvironment{remark}{\begin{rem}\rm}{\end{rem}\rm}

\newcommand{\qqed}{\hspace*{\fill}$\Box$}

\newcommand{\beeq}[1]{\begin{eqnarray}\label{#1}}
\newcommand{\eneq}{\end{eqnarray}}
\newcommand{\cal}{\mathcal}

\newcommand{\ke}{{\cal E}}
\newcommand{\kf}{{\cal F}}
\newcommand{\kh}{{\cal H}}
\newcommand{\kk}{{\cal K}}
\newcommand{\ki}{{\cal I}}
\newcommand{\kj}{{\cal J}}
\newcommand{\kl}{{\cal L}}
\newcommand{\kn}{{\cal N}}

\newcommand{\ko}{{\cal O}}
\newcommand{\kp}{{\cal P}}
\newcommand{\kq}{{\cal Q}}
\newcommand{\kt}{{\cal T}}
\newcommand{\kx}{{\cal X}}

\newcommand{\IC}{{\mathbb C}}

\newcommand{\IP}{{\mathbb P}}
\newcommand{\IQ}{{\mathbb Q}}

\newcommand{\Hom}{{\rm Hom}}
\newcommand{\Ext}{{\rm Ext}}

\newcommand{\End}{{\rm End}}

\renewcommand{\O}       {\mathcal{O}}
\newcommand{\id}{{\rm id}}

\newcommand{\verylongarrow}[1]{\hbox to #1{\rightarrowfill}}

\newcommand\mynote[1]

\newcommand\into{\hookrightarrow}
\begin{document}

\title[$\IP$-twists]{$\IP$-objects and autoequivalences of derived categories}

\author[Daniel Huybrechts]{Daniel Huybrechts}
\address{Mathematisches Institut,
Universit{\"a}t Bonn, Germany}
\email{huybrech@math.uni-bonn.de}

\author[Richard Thomas]{Richard Thomas}
\address{Imperial College, London, UK}
\email{rpwt@ic.ac.uk}

\maketitle


Besides abelian varieties, there are essentially two types of
smooth projective variety with trivial canonical bundle,
Calabi--Yau and holomorphic symplectic manifolds. They are
distinguished, among other things, by their holonomy groups being ${\rm SU}(n)$,
respectively ${\rm Sp}(n)$. The difference between these two
types are in many respect analogous to the difference between
spheres $S^n$ and complex projective spaces
$\IP^n$. This analogy is of course difficult to make precise, but
it comes up in various contexts:

$\bullet$ For a compact manifold $X$ of real dimension $2n$,
respectively $4n$, with holonomy ${\rm SU}(n)$, respectively ${\rm
Sp}(n)$, there exists a ring isomorphism $H^*(X,\ko)\cong H^*({\rm
S}^n,\IC)$, respectively $H^*(X,\ko)\cong H^*(\IP^n,\IC)$.

$\bullet$ The SYZ conjecture predicts in a large complex
structure limit the existence of a lagrangian fibration
$X\twoheadrightarrow S^n$ for any simply connected
Calabi--Yau manifold $X$,
whereas for a holomorphic symplectic manifold this should become
a lagrangian fibration $X\twoheadrightarrow \IP^n$.

$\bullet$ On the symplectic side, Seidel has studied Dehn twists
associated to any lagrangian sphere
contained in a Calabi--Yau manifold, and an analogously defined
symplectomorphism associated to a lagrangian complex projective space (see
\cite[Section 4b]{S1}). In complex dimension 2 (where $S^2\cong\IP^1$), the latter twist is the square of the former; this should be compared to our Proposition \ref{SequalsPinDim1}.

$\bullet$ Spherical twists, studied in detail in \cite{ST}, are
supposed to mirror Dehn twists associated with lagrangian
spheres. These are autoequivalences of the derived category of coherent
sheaves ${\bf D^b}(X)$ of a Calabi--Yau manifold associated to
`spherical' objects $\ke\in{\bf D^b}(X)$. By definition, an
object is spherical if $\Ext^*(\ke,\ke)$ is isomorphic to $H^*(S^n,\IC)$.

\smallskip

This note aims at completing the picture by introducing the
notion of $\IP$-objects and the induced $\IP$-twists, which
are autoequivalences of the derived category ${\bf D^b}(X)$ of a
smooth projective variety $X$ (typically a holomorphic symplectic manifold).

For the low-dimensional case ${\rm SU}(2)={\rm Sp}(1)$ (reflected
by $S^2=\IP^1$) the newly defined $\IP^1$-twist is just the square of the usual spherical twist (see Proposition
\ref{SequalsPinDim1}). In higher dimensions, however,
$\IP^n$-twists describe genuinely new derived equivalences. In
fact, spherical objects ought not exist on manifolds
with holonomy ${\rm Sp}(n\geq2)$.

In higher dimensions one can sometimes nevertheless establish a
relation between $\IP^n$-twists on ${\rm Sp}(n)$-manifolds and
spherical twists on ${\rm SU}(2n+1)$-manifolds. Roughly, a
$\IP^n$-object that does not deform sideways in a one-dimensional
family becomes a spherical object in the ambient
$(2n+1)$-dimensional manifold (see Proposition \ref{Pbecomesph}).
In such a situation Proposition \ref{SversusPtwists} shows that
the spherical twist becomes the $\IP^n$-twist on the special
fibre.

{\bf Acknowledgements:} We would like to thank Paul Seidel for a
useful conversation, Paul Horja, David Ploog and Bal\'azs Szendr\"oi for
comments on the text, and Institut de Math\'ematiques de
Jussieu for support for both authors. The second author is
supported by a Royal Society university research fellowship.

\section{$\IP^n$-objects}

Throughout $X$ will be a smooth projective variety and ${\bf D^b}(X)$ denotes
the bounded derived category of coherent sheaves on $X$. All functors that we use between derived categories are derived functors; we omit ${\bf L}$s and ${\bf R}$s from left and right derived functors.

\begin{definition}
  An object $\ke\in{\bf D^b}(X)$ is called a $\IP^n$-object
  if $\ke\otimes\omega_X\cong\ke$ and $\Ext^*(\ke,\ke)$ is isomorphic
  as a graded ring to $H^*(\IP^n,\IC)$.
\end{definition}

\begin{remark} i) Serre duality shows that if $\ke\in{\bf D^b}(X)$ is
  a $\IP^n$-object then $\dim(X)=2n$.

  ii) The notion of $\IP^n$-objects is designed for holomorphic symplectic
  mani\-folds (i.e.\ hyperk{\"a}hler manifolds). Of course, in this case
  the first condition $\ke\otimes\omega_X\cong\ke$ is automatic.
\end{remark}

\begin{examples} \label{firstegs}
  i) Let $X$ be holomorphic symplectic of dimension $2n$
  and $P:=\IP^n\subset X$. Then
  $\kn_{P/X}\cong \Omega_{P}$ and hence
  $\ke xt^q(\ko_P,\ko_P)\cong\Omega^q_P$. Thus the spectral sequence
  $$E_2^{p,q}=H^p(X,\ke xt^q(\ko_P,\ko_P))\Rightarrow
  \Ext^{p+q}_X(\ko_P,\ko_P)$$ yields an isomorphism
  $\Ext^*_X(\ko_P,\ko_P)\cong H^*(P,\Omega_P^*)=H^*(\IP^n,\IC)$. This is
a ring isomorphism; one roundabout way to see this is given in the remark
following Example \ref{ex}.
  Hence, $\ko_P\in{\bf D^b}(X)$ is a $\IP^n$-object. The same
  arguments also show that $\ko_P(i)\in{\bf D^b}(X)$ is a
  $\IP^n$-object for all $i$.

  ii) Suppose $X$ is an irreducible holomorphic symplectic manifold,
  i.e.\  in addition $H^*(X,\ko_X)\cong H^*(\IP^n,\IC)$. Then any line
  bundle $L$ on $X$ is a $\IP^n$-object. Indeed, $\Ext^*(L,L)\cong
  H^*(X,\ko_X)$.

  iii) Let $X$ be a K3 surface and $C\cong\IP^1\subset X$. Due to i),
  $\ko_C\in{\bf D^b}(X)$ is a $\IP^1$-object. Note that
  $\ko_C$ is also spherical, which reflects $S^2\cong\IP^1$.

  iv) Let $\pi:X\to \IP^n$ be a Lagrangian fibration of an irreducible
  symplectic manifold. If $\ke\in{\bf D^b}(\IP^n)$ is an object with a
graded ring isomorphism
$\bigoplus\Ext^p(\ke,\ke\otimes\Omega_{\IP^n}^q)\cong \bigoplus
H^{p}(\IP^n,\Omega_{\IP^n}^q)$ (in particular, $\ke$ is
exceptional), then
  $\pi^*\ke\in{\bf D^b}(X)$ is a $\IP^n$-object.
  This uses Matsushita's result $R^i\pi_*\ko_X\cong\Omega^i_{\IP^n}$
  \cite{Ma}.  \end{examples}

In dimension $>2$ the notions of a spherical and a $\IP^n$-object are
different. In many examples, $\IP^n$-objects should be thought of as
hyperplane sections of spherical objects. This is made more precise
by the following result.

Suppose ${\cal X}\to C$ is a smooth family over a smooth curve $C$
with distinguished fibre
  $j:X:={\cal X}_0\hookrightarrow{\cal X}$, $0\in C$.
  Let us denote the Kodaira-Spencer class
  of this family by $\kappa({\cal X})\in H^1(X,\kt_X)$ and the Atiyah
  class of an object $\ke\in{\bf D^b}(X)$ by
  $A(\ke)\in\Ext^1(\ke,\ke\otimes\Omega_X^1)$. The product
  $A(\ke)\cdot \kappa({\cal X})\in \Ext^2(\ke,\ke)$ is the
  obstruction to deforming $\ke$ sideways to first order to neighbouring
  fibres in the family \cite{BF,Il}.

\begin{proposition}\label{Pbecomesph}
Suppose $\ke\in{\bf D^b}(X)$ is a $\IP^n$-object such that $0\ne
A(\ke)\cdot \kappa({\cal X})\in \Ext^2(\ke,\ke)$. Then
$j_*(\ke)\in{\bf D^b}({\cal
  X})$ is spherical.
\end{proposition}

\begin{proof}
Firstly, one checks $j_*(\ke)\otimes\omega_{\cal X}\cong
j_*(\ke)$, which follows from $\omega_{\cal X}|_X\cong\omega_X$
and the assumption on $\ke$.

Next, we use the existence of  a distinguished triangle of the form
$$\xymatrix{\ke\otimes\ko_X[1]\ar[r]&j^*j_*\ke\ar[r]&
  \ke\ar[r]&\ke\otimes\ko_X[2]}$$
  with boundary morphism  given by $A(\ke)\cdot \kappa({\cal X}):\ke\to\ke[2]$.
  (For the convenience of the reader we give a short proof of this standard
  result in the appendix.) This yields the long exact sequence
$$\xymatrix{\Ext^k_X(\ke,\ke)\ar[r]&\Ext_{\cal
    X}^k(j_*\ke,j_*\ke)\ar[r]&
  \Ext^{k-1}_X(\ke,\ke)\ar[r]^-\delta&\Ext_X^{k+1}(\ke,\ke),}$$
where we use adjunction $\Ext^k_{\cal
  X}(j_*\ke,j_*\ke)=\Ext_X^k(j^*j_*\ke,\ke)$.

The boundary morphism $\delta$ is given by cup-product with
$A(\ke)\cdot \kappa({\cal X}):\ke\to\ke[2]$ considered as an element in
$\Ext^2(\ke,\ke)$.  By assumption this may be taken to be the degree 2
generator of $\Ext^*(\ke,\ke)$. Therefore $\delta$
is an isomorphism for $1\le k\le 2n-1$, yielding

$$\Ext^k(j_*\ke,j_*\ke)=\left\{\begin{array}{ll}
  \IC&k=0,2n+1,\\0&{\rm otherwise.}\end{array}\right.$$
Thus, $j_*\ke$    is indeed a spherical object in ${\bf D^b}({\cal
X})$ of the $(2n+1)$-dimensional variety ${\cal X}$.
\end{proof}

\begin{example}\label{ex}
  The typical example for the situation considered above is the
  twistor  space of a hyperk{\"a}hler manifold. The twistor space is,
  unfortunately,  never a projective manifold (not even K{\"a}hler),
  but the notion of a spherical object makes sense also in the
  analytic category.

  For the $\IP^n$-object $\ko_P$ given by a projective space $P\subset
  X$ one can in fact find a projective $\kx$  for which $\IP^n\cong P\subset
  X$ does not deform to $P_t\subset {\cal X}_t$ even to first order
  (see \cite{HDG}).

  Note that there are $\IP^n$-objects for which such a family
  does not exist. E.g.\ the $\IP^n$-object $\ko_X$ always deforms
  in families. On the other hand, any non-trivial line bundle
admits such a family.

Presumably this might be fixed by deforming $X$ in noncommutative directions, i.e.\ allowing ${\cal X}$ to be a noncommutative variety, providing
other examples of the same phenomenon, namely that a
$\IP^n$-object becomes a spherical object in an `ambient' derived
category (the analogue of ${\bf D^b}({\cal X})$).
\end{example}

We remark that the proof of Proposition \ref{Pbecomesph} shows
that if $\Ext^*(\ke,\ke)$ is isomorphic to $H^*(\IP^N,\IC)$ as a
vector space but \emph{not} as a ring, then $j_*\ke$ is
\emph{not} a spherical object. In the case of $\ko_P$ given by a
projective space $\IP^N\cong P\subset X$ with normal bundle
$\Omega_P$ as in Example \ref{firstegs} i), fix an $\kx$ for which
$P\subset X$ does not deform to first order. It is then standard
that the normal bundle of $P\subset\kx$ is isomorphic to
$\ko_{\IP^n}(-1)^{\oplus(n+1)}$, from which it is easy to see it
is spherical (see e.g.\ \cite[Ch.\ 8]{HFM}). Thus the ring
structure on $\Ext^*(\ko_P,\ko_P)$ is indeed that of
$H^*(\IP^N,\IC)$.

\section{$\IP^n$-twists}

We shall try to imitate the construction of the spherical twist
$T_\ke$ associated to any spherical object (see \cite{ST}) and
define a $\IP^n$-twist for any $\IP^n$-object. This is done in
two steps. We first describe the Fourier--Mukai kernel and then
show that the induced Fourier--Mukai transform is an equivalence.
The second step is straightforward, whereas the description of the
Fourier--Mukai kernel itself is interesting in as much as it uses
a double cone construction.

Suppose $\ke\in{\bf D^b}(X)$ is a $\IP^n$-object. A
generator
$\bar{h}\in\Ext^2(\ke,\ke)$
will be viewed as a morphism $h:\ke[-2]\to\ke$.
The ring $\Ext^*(\ke,\ke)$ is then isomorphic to $\IC\left[\bar{h}\right]/(\bar
h^{n+1})$.
The image of $\bar{h}$ under the natural isomorphism $\Ext^2(\ke,\ke)\cong
\Ext^2(\ke\dual,\ke\dual)$ will be denoted $\bar{h}\dual$, which
represents a morphism $h\dual:\ke\dual[-2]\to\ke\dual$.

Then introduce
$H:=h\dual\boxtimes{\rm id}-{\rm id}\boxtimes h$  on $X\times X$ which is thus a morphism
$$\xymatrix{H:\left(\ke\dual\boxtimes\ke\right)[-2]\ar[r]&\ke\dual\boxtimes\ke.}$$
The cone $\kh:={\rm C}(H)$ of this morphism
fits in a distinguished
triangle
$$\xymatrix{\left(\ke\dual\boxtimes\ke\right)[-2]\ar[r]^-H&\ke\dual\boxtimes\ke\ar[r]&\kh\ar[r]&\left(\ke\dual\boxtimes\ke\right)[-1].}$$

Recall that the kernel of the spherical twist associated to a
spherical object is by definition the cone of the
trace morphism ${\rm tr}:\ke\dual\boxtimes\ke\to\ko_\Delta$, which
is actually the composition of the restriction to the diagonal
$\Delta\subset X\times X$ and the
trace map on the diagonal.

In order to define the kernel of the $\IP^n$-twist we have to go one
step further.

%

\begin{lemma}
  The natural trace map ${\rm tr}:\ke\dual\boxtimes\ke\to\ko_\Delta$
  factorizes   uniquely over the cone $ \kh$, i.e.\
  there exists a unique morphism ${\rm t}$ that  makes the following diagram commutative
  $$\xymatrix{\ke\dual\boxtimes\ke\ar[dr]_{\rm
      tr}\ar[r]&\kh\ar[d]^{{\rm t}}\\
  &~~~\ko_\Delta.}$$
\end{lemma}

\begin{proof}
  Apply $\Hom(~~~,\ko_\Delta)$ to the distinguished triangle defining $\kh$.
  Use $\Ext^i_{X\times X}(\ke\dual\boxtimes\ke,\ko_\Delta)
  \cong \Ext^i_X(\ke,\ke)$ and the definition of $H$, to show that
  the boundary maps   $\Ext^i(\ke,\ke)\to \Ext^{i+2}(\ke,\ke)$
  become $\bar h-\bar h=0$. Hence $\Hom(\kh,\ko_\Delta)\to
\Hom(\ke\dual\boxtimes\ke,\ko_\Delta)$
is an isomorphism, giving the unique lift ${\rm t}$ of the trace map.
\end{proof}

\begin{definition}
  To any $\IP^n$-object $\ke\in{\bf D^b}(X)$ one associates the cone
  $$\kq_\ke:=
  {\rm C}(\!\xymatrix{\kh\ar[r]^-{\rm t}&\ko_\Delta}\!)\in{\bf D^b}(X\times
  X).$$
\end{definition}

While we used the assumption that $\ke$ is a $\IP^n$-object to simplify 
the above construction, an alternative construction, using locally free 
and \v Cech resolutions, shows that in fact t always exists.

\begin{definition}
Let $\ke\in{\bf D^b}(X)$ be a $\IP^n$-object and let
$\kq_\ke\in{\bf D^b}(X\times X)$ be the object associated to it
by the above construction. The $\IP^n$-twist $P_\ke$ induced by a
$\IP^n$-object $\ke\in{\bf D^b}(X)$ is the Fourier--Mukai
transform
$$\xymatrix{P_\ke:=\Phi_{\kq_\ke}:{\bf D^b}(X)\ar[r]&{\bf D^b}(X)}$$
with kernel $\kq_\ke$.
\end{definition}

\begin{rem}
Note that the induced actions $P_\ke^K:K(X)\to K(X)$ and
$P_\ke^H:H^*(X,\IQ)\to H^*(X,\IQ)$ are both the identity. This follows
from the observation that
$[\kh]=[\ke\dual\boxtimes\ke]+[\ke\dual\boxtimes\ke[-1]]=0$ in
K-theory. Later, in Remark \ref{defid}, we will see that in most cases $P_\ke$
can even be deformed to the identity on a deformation of $X$.
\end{rem}

It is often useful to have a cone description of the image of an
object under the $\IP^n$-twist. It is not difficult to see that
for any $\kf\in {\bf D^b}(X)$ the image $P_\ke(\kf)$ is isomorphic
to the double cone
$${\rm C}\left(\!\!\xymatrix{{\rm C}\big(\Ext^{*-2}(\ke,\kf)\otimes
\ke\ar[rr]^{\ \bar h\dual\!\cdot\id-\id\cdot h}&&
\Ext^*(\ke,\kf)\otimes\ke\big)\ar[r]&\kf}\!\right).$$

Let us spell this out in  the case
of the $\IP^n$-object itself and objects that are orthogonal to
it.

\begin{lemma}
  Let $\ke\in{\bf D^b}(X)$ be a $\IP^n$-object. Then

  {\rm i)} $P_\ke(\ke)\cong \ke[-2n]$, and

{\rm ii)} $P_\ke(\kf)\cong\kf$ for any $\kf\in\ke^\perp:=\{\mathcal G\colon\Ext^*(\ke,\mathcal G)=0\}$.
\end{lemma}

\begin{proof}
The second assertion is trivial, for $\Phi_\kh(\kf)\cong0$ for any
$\kf\in\ke^\perp$ and, therefore,
$P_\ke(\kf)\cong\Phi_{\ko_\Delta}(\kf)$.

For the first one we use
$\Phi_{\ke\dual\boxtimes\ke}(\ke)=\bigoplus \bar{h}^i\cdot\ke[-2i]$
to compute  $\Phi_{\kh}(\ke)$ as the cone of the morphism
$$\xymatrix{\id\cdot\ke[-2]\ar@{}[d]|\oplus\ar[rr]^{-h}\ar[drr]^{\bar{h}\cdot}&&\id\cdot\ke\ar@{}[d]|\oplus\\
  \bar{h}\cdot\ke[-4]\ar@{}[d]|\oplus\ar[rr]^{-h}\ar[drr]^{\bar{h}\cdot}&&\bar{h}\cdot\ke[-2]\ar@{}[d]|\oplus\\
  \vdots\ar@{}[d]|\oplus&&\vdots\ar@{}[d]|\oplus\\
  \bar{h}^n\cdot\ke[-2n-2]\ar[rr]^{-h}&&\bar{h}^n\cdot\ke[-2n].}$$
  The evaluation map to $\ke$ is the obvious one, taking the $r$th term in
  the right hand column to $\ke$ by the map $h^r$, for all $0\le r\le n$.

  $C\big(\ke\stackrel{\id\,}{\to}\ke\big)$ maps into the cone on the evaluation map
  in the obvious way (with the first factor mapping isomorphically to the
  top right hand corner of the above diagram). Taking the cone on this shows
  that $P_\ke(\ke)$ is quasi-isomorphic to the cone on the above diagram with the
  top right hand corner removed. Repeating this procedure with the subcone $C\big(\id\cdot\ke[-2]\stackrel{\bar h\cdot}{\to}\bar h\cdot\ke[-2]\big)$ shows we can
further remove the top left hand corner and the next right hand term. Iterating
leaves us with just the bottom left hand term $\bar{h}^n\cdot\ke[-2n-2]$
which, due to its position, makes $P_\ke(\ke)\cong \ke[-2n]$.
\end{proof}

\begin{proposition}
  For any $\IP^n$-object $\ke\in{\bf D^b}(X)$ the associated
  $\IP^n$-twist is an autoequivalence
  $$\xymatrix{P_\ke:{\bf D^b}(X)\ar[r]^-\sim&{\bf D^b}(X).}$$
\end{proposition}

\begin{proof}
  We follow Ploog's simplified proof (see \cite{Ploog}) of the analogous
  result for spherical twists in \cite{ST}.

  Let $\ke\in{\bf D^b}(X)$ be any object. Then
  $\Omega:=\{\ke\}\cup\ke^\perp$ is a spanning class
(cf.\ \cite[Ch.8]{HFM}). On this spanning
  class the $\IP^n$-twist acts by shifting $[-2n]$ on $\ke$ and as the
  identity  on the rest. This immediately shows that $P_\ke$ is fully
  faithful.

  In order to show that $P_\ke$ is an equivalence use the assumption
  $\ke\otimes\omega_X\cong\ke$ which ensures that $\kq_\ke\otimes
  \pi_1^*\omega_X\cong\kq_\ke\otimes\pi_2^*\omega_X$, where $\pi_i$ is the
projection of $X\times X$ onto its $i$th factor. Hence the left and right
  adjoint of $P_\ke$ coincide, which suffices to conclude
(cf.\ \cite{B,HFM}).
\end{proof}

Let us compare $\IP^n$-twists and spherical twists. We study the
situation of Proposition \ref{Pbecomesph}. So, let ${\cal X}\to C$ be
a smooth family over a smooth curve $C$ with distinguished fibre
  $j:X:={\cal X}_0\hookrightarrow{\cal X}$, $0\in C$.

\begin{proposition}\label{SversusPtwists}
Suppose $\ke\in{\bf D^b}(X)$ is a $\IP^n$-object with
$A(\ke)\cdot\kappa({\cal X})\ne0$. Then $j_*$ intertwines the $\IP^n$-twist $P_\ke$ and the spherical twist $T_{j_*\ke}$, i.e.\ one has the
following commutative diagram
$$\xymatrix{{\bf D^b}(X)\ar[d]_{P_\ke}\ar[rr]^{j_*}&&{\bf D^b}({\cal
    X})\ar[d]^{T_{j_*\ke}}\\
  {\bf D^b}(X)\ar[rr]^{j_*}&&{\bf D^b}({\cal
    X}).}$$
\end{proposition}

\begin{proof}
  This is an application of Chen's lemma (see \cite{Chen} or \cite[Ch.\ 11]{HFM}). One simply has to show
  that there exists an object $\kl$ on ${\cal X}\times_C{\cal X}$ with
  $$\ell_*\kl\cong \kp_{j_*\ke}~~~~~~{\rm and}~~~~~~
  f^*\kl\cong\kq_{\ke}.$$
  Here $\kp_{j_*\ke}={\rm C}(j_*(\ke)\dual\boxtimes
  j_*\ke\to\ko_\Delta)$ is the kernel of the spherical twist
  $T_{j_*\ke}$, and notations for the relevant morphisms are fixed as follows
  $$\xymatrix{X\ar[r]^j\ar[d]_{\iota_0}&{\cal
  X}\ar[d]_k\ar[dr]^\iota&\\
  X\times X\ar[r]_f&{\cal X}\times_C{\cal
      X}\ar[r]_\ell&{\cal X}\times{\cal X\,,}}$$
 where $\iota_0$, $\iota$, and $k$ are the diagonal embeddings.

  We shall define $\kl$ as the cone ${\rm C}(\xi)$ of a morphism
  $$\xymatrix{\xi:f_*(\ke\dual[-1]\boxtimes\ke)\ar[r]& k_*\ko_{\cal X}}.$$

  The morphism $\xi$ itself is constructed as a composition as follows. Consider first
  the trace  map
  ${\rm tr}:\ke\dual[-1]\boxtimes\ke\to {\iota_0}_*\ko_X[-1]$ and take its image under $f_*$.
  Then use the short exact sequence $0\to\ko_{\cal X}\to\ko_{\cal X}(X)\to
  j_*\ko_X\to0$ and the
  induced  boundary map $j_*\ko_X[-1]\to\ko_{\cal X}$.
   Using $k\circ j=f\circ \iota_0$, the image of the latter under
   $k_*$ can be composed with $f_*({\rm tr})$. This yields
\begin{equation}\label{xi}
\xymatrix{\xi:f_*(\ke\dual[-1]\boxtimes\ke)
    \ar[r]& f_*{\iota_0}_*\ko_X[-1]\cong k_*j_*\ko_X[-1]\ar[r]&
  k_*\ko_{\cal X}.}
\end{equation}

Its cone $\kl:={\rm C}(\xi)$ can be alternatively described as
\begin{eqnarray*}
\ell_*\kl&\cong&\ell_*{\rm C}(\xymatrix {
f_*(\ke\dual[-1]\boxtimes\ke)\ar[r]^-{\xi}&
k_*\ko_{\cal X}})\\
&\cong&{\rm
  C}(\xymatrix{\ell_*f_*(\ke\dual[-1]\boxtimes\ke)\ar[r]^-{\ell_*\xi}&\ell_*k_*\ko_{\cal X}})\\
  &\cong&{\rm C}(\xymatrix{(j_*\ke)\dual\boxtimes j_*\ke\ar[r]&\iota_*\ko_{\cal
  X}}),
\end{eqnarray*}
where we use duality for the isomorphism
$j_*(\ke\dual[-1])\cong(j_*\ke)\dual$. Also observe that the
functorial properties of duality imply that indeed $\ell_*\xi={\rm
tr}$ and hence $\ell_*\kl\cong \kp_{j_*\ke}$.

To show $f^*\kl\cong\kq_\ke$ one first observes that
$f^*k_*\ko_{\cal X}\cong{\iota_0}_*\ko_X$, because the
intersection of $f(X\times X)$ and $k({\cal X})$ (inside ${\cal
X}\times_C{\cal X}$) is transversal. Next we use the existence of
the distinguished triangle, for an object on  the divisor $f:X\times
X\hookrightarrow{\cal X}\times_C{\cal X}$,
$$\xymatrix{
(\ke\dual\boxtimes\ke)[-2]\ar[r]^-\delta&\ke\dual\boxtimes\ke\ar[r]&f^*f_*(\ke\dual[-1]\boxtimes\ke)\ar[r]&
(\ke\dual\boxtimes\ke)[-1]}\!,$$
with the boundary map $\delta$ given by the cup-product with the
obstruction class $A(\ke\dual\boxtimes\ke)\cdot \kappa({\cal
X}\times_C{\cal X})\in\Ext^2_{X\times
X}(\ke\dual\boxtimes\ke,\ke\dual\boxtimes\ke)$ (see the appendix).

As we may assume that $A(\ke)\cdot\kappa({\cal X})=\bar{h}$ and passing
from a bundle to its dual changes the Atiyah class by a sign, one
finds $A(\ke\dual\boxtimes\ke)\cdot \kappa({\cal X}\times_C{\cal
X})={\rm id}\boxtimes \bar{h}-\bar{h}\dual\boxtimes{\rm id}$. Hence
$f^*f_*(\ke\dual[-1]\boxtimes\ke)\cong \kh$.

The last thing one has to check before concluding
$f^*\kl\cong\kq_\ke$ is the commutativity of
$$\xymatrix{\ke\dual\boxtimes\ke\ar[r]\ar[dr]_{\rm
tr}&f^*f_*(\ke\dual[-1]\boxtimes\ke)\ar[d]^{f^*\xi}\\
&{\iota_0}_*\ko_X,}$$
i.e. that the following diagram, whose vertical part is
(\ref{xi}), commutes
\begin{equation}\label{rectangle}
\xymatrix{\ke\dual\boxtimes\ke\ar[r]\ar[d]_{\rm
tr}&f^*f_*(\ke\dual[-1]\boxtimes\ke)\ar[d]^{f^*f_*{\rm tr}}
\ar@/^5pc/[dd]^{f^*\xi}
\\
{\iota_0}_*\ko_X\ar[r]\ar[dr]_-\cong&f^*f_*{\iota_0}_*\ko_X[-1]\ar[d]\\
&f^*k_*\ko_{\cal X}.}
\end{equation}
Commutativity of the rectangle follows from the functoriality of
the distinguished triangle constructed in Proposition \ref{AtKod}
in the appendix, i.e. from the commutativity of
(\ref{commutative}) applied to $\ke_1=\ke\dual\boxtimes\ke,
\,\ke_2=\iota_{0*}\ko_X\in{\bf D^b}(X\times X)$ and the divisor
$f:X\times X\to\kx\times_C\kx$.

To show commutativity of the triangle, we apply $f^*$ to
$0\to k_*\O_\kx\to k_*\O_\kx\to f_*\iota_{0*}\O_X\to0$. This gives an exact triangle
$$
\iota_{0*}\O_X\to f^*f_*\iota_{0*}\O_X\stackrel{e\,}{\to}
\iota_{0*}\O_X[1]\to\iota_{0*}\O_X[1],
$$
whose boundary morphism is zero (since
$L_0f^*f_*\iota_{0*}\O_X\cong\iota_{0*}\O_X$, for instance). Thus
$f^*f_*\iota_{0*}\O_X\cong\iota_{0*}\O_X\oplus\iota_{0*}\O_X[1]$,
to which the only morphisms from $\iota_{0*}\O_X[1]$ are
multiplies of $(0,\id)$, since $\iota_{0*}\O_X$ is a simple
\emph{sheaf}. Therefore in our exact triangle of Proposition
\ref{AtKod}
$$
\iota_{0*}\O_X[1]\to f^*f_*\iota_{0*}\O_X\to\iota_{0*}\O_X\to\iota_{0*}\O_X[2].
$$
the first arrow must be a nonzero multiple of $(0,\id)$, and so the above
morphism $e:f^*f_*\iota_{0*}\O_X\to\iota_{0*}\O_X[1]$ splits the exact triangle.

If we show that this morphism $e$ is the vertical arrow of the triangle (\ref{rectangle})
then we have shown that the composition of the horizontal and vertical morphisms in the triangle is the natural identification between $\iota_{0*}\O_X$ and $f^*k_*\ko_{\cal X}$, i.e.\ the triangle is commutative.

So finally we must check in the definition of $\xi$ that
$$f^*k_*(j_*\ko_X[-1]\to\ko_\kx)=f^*\left(k_*\ko_\kx\otimes
\left[f_*\ko_{X\times X}[-1]\to\ko_{\kx\times_C\kx}\right]\right).$$ This follows
from the fact that the natural map between the short exact
sequences
$$\xymatrix{0\ar[r]&\ko_{\kx\times_C\kx}\ar[r]\ar[d]&\ko_{\kx\times_C\kx}\ar[d]\ar[r]&f_*\ko_{X\times X}\ar[d]\ar[r]&0\\
0\ar[r]&k_*\ko_\kx\ar[r]&k_*\ko_\kx\ar[r]&k_*j_*\ko_X\ar[r]&0}$$
factorizes over the tensor product of the first one with
$k_*\ko_\kx$, which stays exact due to the transversality of the
intersection of $f(X\times X)$ and $k(\kx)$.
\end{proof}

\begin{rem}\label{defid}
The restriction of the spherical kernel of $j_*\ke$ to a fibre
${\cal X}_t$ with $t\ne0$ is isomorphic to the diagonal and the
induced equivalence ${\bf D^b}({\cal X}_t)\cong{\bf D^b}({\cal
X}_t)$ is the identity. In this sense, a projective twist
associated to a $\IP^n$-object $\ke$ with $A(\ke)\cdot
\kappa({\cal X})\ne0$ can be deformed to the identity.

This is the derived version of the fact that any birational
correspondence between holomorphic symplectic varieties can be
deformed to an isomorphism (see \cite{HDG}). It is also mirror to
what Seidel calls `fragility' of his Dehn twists about lagrangian
$\IP^n$ submanifolds, at least in the case $n=1$ (see \cite{S2}).
\end{rem}

Let us conclude with a discussion of the  two-dimensional
situation, where spherical and $\IP^1$-twist are related more
directly.

\begin{proposition}\label{SequalsPinDim1}
  Let $\ke\in{\bf D^b}(X)$ be a $\IP^1$-object (thus $\dim(X)=2$).
  Then $$T^2_{\ke}\cong P_\ke.$$
\end{proposition}

\begin{proof}
  In order to compare the Fourier-Mukai kernel
  \begin{equation}\label{FMP1}{\rm C}\left({\rm C}
      (\ke\dual\boxtimes\ke[-2]\stackrel{H\,}{\to}\ke\dual\boxtimes\ke))
\to\ko_\Delta\right)
      \end{equation}
 of the $\IP^1$-twist $P_\ke$ with the kernel $\kk$ of the square $T^2_\ke$,
 we will compute the latter explicitly by the standard product formula.
 Let ${\rm C}:={\rm C}(\ke\dual\boxtimes\ke\to\ko_\Delta)$. Then
 $\kk\cong\pi_{13*}(\pi_{12}^*{\rm C}\otimes\pi_{23}^*{\rm C})$,
    with $\pi_{ij}:X\times X\times X\to X\times X$ denoting the usual
    projections.

    The tensor product involves terms in degree $-2,-1$, and $0$.
    In degree $-2$ this is $\kl_2:=\pi_{12}^*(\ke\dual\boxtimes\ke)\otimes
    \pi_{23}^*(\ke\dual\boxtimes\ke)$, in degree $-1$ one
    finds $\kl_1:=(\pi_{12}^*(\ke\dual\boxtimes\ke)\otimes\pi_{23}^*\ko_\Delta)
    \oplus(\pi_{12}^*\ko_\Delta\otimes\pi_{23}^*(\ke\dual\boxtimes\ke))$
    and in degree $0$ simply
    $\kl_0:=\pi_{12}^*\ko_\Delta\otimes\pi_{23}^*\ko_\Delta$.

    Clearly,
    $\pi_{13*}\kl_2\cong(\ke\dual\boxtimes\ke)\otimes\End^*(\ke)$,
    $\pi_{13*}\kl_1\cong(\ke\dual\boxtimes\ke)\oplus(\ke\dual\boxtimes
    \ke)$,
    and $\pi_{13*}\kl_0\cong\ko_\Delta$. Moreover,
    $\pi_{13*}\kl_1\to\pi_{13*}\kl_0$ is ${\rm tr}\oplus{\rm tr}$ and
    $\pi_{13*}\kl_2\cong(\ke\dual\boxtimes\ke)\oplus(\ke\dual\boxtimes\ke)[-2]\to\pi_{13*}\kl_1$
    is the diagonal on the first summand and $h\dual\boxtimes{\rm
    id}\oplus{\rm id}\boxtimes h$ on the second.

    To conclude, embed the complex
$\ke\dual\boxtimes\ke\stackrel{\id\,}{\to}\ke\dual\boxtimes\ke$ into
$\pi_{13*}\kl_2\to\pi_{13*}\kl_1$ via
$(1,0):\ke\dual\boxtimes\ke\to\pi_{13*}\kl_2=(\ke\dual\boxtimes\ke)\oplus
(\ke\dual\boxtimes\ke[-2])$ and
$(1,1):\ke\dual\boxtimes\ke\to\pi_{13*}\kl_1$ respectively. The cokernel of this map is identified with
$\ke\dual\boxtimes\ke[-2]\stackrel{H\,}{\to}\ke\dual\boxtimes\ke$
by means of the second projection
$\pi_{13*}\kl_2\to\ke\dual\boxtimes\ke[-2]$ and
$(1,-1):\pi_{13*}\kl_1=(\ke\dual\boxtimes\ke)\oplus(\ke\dual\boxtimes\ke)\to\ke\dual\boxtimes\ke$.
Thus $\kk$ is isomorphic to the kernel (\ref{FMP1}).
\end{proof}

\section{Appendix}

Let $\kx\to C$ be a smooth projective morphism  over a smooth
curve  $C$ with parameter $t$. The central fibre  will be called
$X=\kx_0$ and its inclusion  $j: X\into\kx$.

The family $\kx$ viewed as a deformation of $X$ induces the
Kodaira-Spencer class $\kappa(\kx)\in H^1(X,\kt_X)$, which is by
definition the extension class of the normal bundle sequence
\begin{equation}\label{KSpre}
\xymatrix{0\ar[r]&\kt_X\ar[r]&\kt_\kx|_X\ar[r]&\ko_X\ar[r]&0.}
\end{equation}
(Of course multiplication by $t$ induces a trivialization of the
normal bundle: $\ko_X\cong\ko_X(X)$.) The sequence can be dualized
to yield
\begin{equation}\label{KS}
\xymatrix{0\ar[r]&\ko_X\ar[r]&\Omega_\kx|_X\ar[r]&\Omega_X\ar[r]&0,}
\end{equation}
and the Kodaira-Spencer class will be viewed as its boundary
morphism $\kappa(\kx):\Omega_X\to\ko_X[1]$.

For $\ke\in {\bf D^b}(X)$ we denote by $J(\ke)$ its first jet
space, i.e.
$$
J(\ke)=\pi_{2*}\big(\pi_1^*\ke\otimes\ko_{2\Delta}\big),
$$
where $\Delta\subset X\times X$ is the diagonal, $2\Delta$ is its
double: $\ki_{2\Delta}:=\ki_\Delta^2$, and $\pi_i$ is the
projection onto the $i$th factor $X$.

As $\ko_{2\Delta}$ sits in the short exact sequence
$$\xymatrix{0\ar[r]&\Omega_\Delta\ar[r]&\ko_{2\Delta}\ar[r]&\ko_\Delta\ar[r]&0,}$$
the jet space $J(\ke)$ sits in a distinguished triangle of the
form
\begin{equation}\label{Jet}
\xymatrix{\ke\otimes\Omega_X\ar[r]&J(\ke)\ar[r]&\ke\ar[r]&\ke\otimes\Omega_X[1].}
\end{equation}
The extension class, i.e.\ the boundary morphism
$\ke\to\ke\otimes\Omega_X[1]$, is by definition the Atiyah class
$A(\ke)\in \Ext^1(\ke,\ke\otimes\Omega_X)$.

The product $A(\ke)\cdot \kappa({\cal X})\in
\Ext^2(\ke,\ke)=\Hom(\ke[-1],\ke[1])$ can be
described as the composition of
$A(\ke)[-1]:\ke[-1]\to\ke\otimes\Omega_X$ with ${\rm
id}_\ke\otimes\kappa(\kx):\ke\otimes\Omega_X\to\ke[1]$. In
particular, there exists a distinguished triangle
\begin{equation}\label{AKJet}
\xymatrix{\ke[1]\ar[r]&{\rm C}\left(\ke\otimes\Omega_\kx|_X\to
J(\ke)\right)\ar[r]&\ke\ar[r]^-{A\cdot\kappa}&\ke[2]}.
\end{equation}

\begin{proposition}\label{AtKod}
Let $\ke[-1]\to\ke[1]$ be the morphism given by $A(\ke)\cdot \kappa({\cal X})$ as above. Then there exists a functorial (in $\ke$)
isomorphism ${\rm C}(\ke[-1]\to\ke[1])\cong j^*j_*\ke$.
\end{proposition}

Note that in the general situation of a divisor $j:X\into \kx$
there always exists a distinguished triangle of the form (see
e.g.\ \cite[Ch.\ 11]{HFM})
\begin{equation}\label{classic}
\xymatrix{\ke\otimes\ko_X(-X)[1]\ar[r]&j^*j_*\ke\ar[r]&\ke\ar[r]&\ke\otimes\ko_X(-X)[2]}.\end{equation}
and the proposition asserts that in the special case of a family
$\kx\to C$ there exists such a triangle with the boundary
morphism given by $A(\ke)\cdot \kappa({\cal X})$. (Of course the assertion
holds true also in the general situation, with $\kappa(\kx)$
defined appropriately as a class in
$\Ext^1(\Omega_X,\ko_X(-X))$, but we won't need this.)

In the following we shall denote by $\ki$ the ideal sheaf of the
diagonal $\iota:\Delta\hookrightarrow X\times X$ and by $\kj$ the
ideal sheaf of $\Delta$ as a subvariety of $X\times\kx$ via the
closed embedding $i:={\rm id}\times j:X\times X\hookrightarrow
X\times\kx$. Then $\iota_*\Omega_X\cong\ki/\ki^2$,
$\iota_*(\Omega_\kx|_X)\cong i^*(\kj/\kj^2)$, and the conormal
bundle sequence of $X\hookrightarrow \kx$ is induced by the
natural map $i^*(\kj/\kj^2)\to \ki/\ki^2$. The latter can be
composed with $\ki/\ki^2\to\ko_{X\times X}/\ki^2=\ko_{2\Delta}$
to yield $\eta:i^*(\kj/\kj^2)\to\ko_{2\Delta}$. The cone of this
morphism $\iota_*(\Omega_\kx|_X)\stackrel{\eta\,}{\to}\ko_{2\Delta}$
is described by the following lemma.

\begin{lemma}
There is a natural isomorphism ${\rm C}(\eta)\cong
i^*i_*(\iota_*\ko_\Delta)$.
\end{lemma}

\begin{proof}
Pulling-back the short exact sequence
$$\xymatrix{0\ar[r]&\kj\ar[r]&\ko_{X\times\kx}\ar[r]&i_*\iota_*\ko_\Delta\ar[r]&0}$$
via $i$ yields a distinguished triangle
$$\xymatrix{i^*\kj\ar[r]&\ko_{X\times
X}\ar[r]&i^*i_*(\iota_*\ko_\Delta)\ar[r]&i^*\kj[1].}$$ (Note that
$i^*\kj$ need not be derived; it is the normal pullback. This follows, for instance, from the fact that $i^*i_*(\iota_*\ko_\Delta)$ has
cohomology only in degrees $0$ and $-1$, by (\ref{classic}).)

Then use the diagram
$$\xymatrix{\kk\ar[d]\ar[r]&\ki^2\ar[d]&\\
i^*\kj\ar[d]\ar[r]&\ko_{X\times X}\ar[r]\ar[d]&i^*i_*(\iota_*\ko_\Delta)\\
i^*(\kj/\kj^2)\ar[r]&\ko_{2\Delta}&}$$ and the fact, which can be
easily verified by a local calculation, that the natural map
$$\kk:=\ker(i^*\kj\to i^*(\kj/\kj^2))\to\ki^2$$
is an isomorphism. Hence, $$i^*i_*(\iota_*\ko_\Delta)\cong{\rm
C}(i^*\kj\to\ko_{X\times X})\cong {\rm
C}(i^*(\kj/\kj^2)\to\ko_{2\Delta}).$$
\end{proof}

{\it Proof of proposition.} Denote by $p_i$  the projection from
$X\times\kx$ onto the $i$th factor. Then use $p_2\circ
i=j\circ\pi_2$ to conclude
$$j^*j_*\ke\cong j^*j_*\pi_{2*}(\pi_1^*\ke\otimes \iota_*\ko_\Delta)\cong
j^*p_{2*}i_* (\pi_1^*\ke\otimes \iota_*\O_\Delta).$$

Use basechange  for $p_2\circ i=j\circ \pi_2$
(see \cite[Lemma 1.3]{BO}) to identify this
with
$$
\pi_{2*}i^*i_*(\pi_1^*\ke\otimes \iota_*\O_\Delta)\cong
\pi_{2*}i^*(p_1^*\ke\otimes i_*(\iota_*\O_\Delta)).
$$
Now $p_1\circ i=\pi_1$ and tensor product commutes with pullback,
so this is
$$
\pi_{2*}\big((\pi_1^*\ke)\otimes i^*i_*(\iota_*\O_\Delta)\big).
$$
So by the lemma,
\begin{eqnarray*} j^*j_*\ke&\cong& {\rm C}\left(
\pi_{2*}(\pi_1^*\ke\otimes i^*(\kj/\kj^2))\to\pi_{2*}(
\pi_1^*\ke\otimes\O_{2\Delta})\right)\\
&\cong& {\rm C}\left(\ke\otimes\Omega_\kx|_X\to J(\ke)\right).
\end{eqnarray*}
Functoriality is clear from the functoriality of pushforwards and pullbacks,
which is all that we have used. That is a morphism $\varphi:\ke_1\to\ke_2$ in ${\bf D^b}(X)$ induces a commutative diagram
\begin{equation}\label{commutative}
\xymatrix{\ke_1[1]\ar[d]_{\varphi[1]}\ar[r]& j^*j_*\ke_1\ar[d]_{j^*j_*\varphi}\ar[r]&\ke_1\ar[d]^\varphi
\ar[rr]^-{A(\ke_1)\cdot\kappa(\kx)}&&\ke_1[2]\ar[d]^{\varphi[2]}\\
\ke_2[1]\ar[r]& j^*j_*\ke_2\ar[r]&\ke_2\ar[rr]^-{A(\ke_2)\cdot\kappa(\kx)}&&\ke_2[2].}
\end{equation}
\qqed

\bigskip

 

\end{document}